\documentclass{article}
\RequirePackage[text={140mm,191mm},centering,paperwidth=8.5in,paperheight=11.0in]{geometry}
\usepackage[defaultlines=2,all]{nowidow}
\usepackage{wrapfig}

\usepackage{graphicx}        
\usepackage[bottom]{footmisc}

\usepackage{amsmath}
\usepackage{amsthm}
\usepackage{amsfonts}
\usepackage{amssymb}
\usepackage{mathrsfs}
\usepackage{url}
\usepackage{tikz}
\usetikzlibrary{calc}
\usetikzlibrary{decorations.markings}
\usepackage{alltt}
\usepackage{enumitem}
\usepackage{authblk}
\RequirePackage{etoolbox}
\makeatletter
\patchcmd{\@maketitle}{center}{flushleft}{}{}
\patchcmd{\@maketitle}{center}{flushleft}{}{}
\patchcmd{\@maketitle}{\LARGE}{\LARGE\bf\boldmath}{}{}

\def\maketitle{{%
  \renewenvironment{tabular}[2][]
    {\begin{flushleft}}
    {\end{flushleft}}
  \AB@maketitle}}
\makeatother

\RequirePackage{hyperref}
\hypersetup{
    colorlinks=true,
    linkcolor=black,
    filecolor=black,      
    urlcolor=blue,
    citecolor=black,
} 
\newcommand*{\email}[1]{\hspace*{0.5em}{\rm\href{mailto:#1}{#1}}}

\theoremstyle{thmstyleone}%
\newtheorem{theorem}{Theorem}

\theoremstyle{thmstyletwo}%

\theoremstyle{thmstylethree}%

\makeindex             

\begin{document}

\title{Small transitive homogeneous $3$-$(v,\{4,6\},1)$ designs}
\author[1]{Michael Epstein}
\author[2]{Donald L. Kreher}
\author[3]{Spyros S. Magliveras}
\affil[1]{Colorado State University, Fort Collins, Colorado,  80523 USA\newline
\email{Michael.Epstein@colostate.edu}
}
\affil[2]{Michigan Technological University, Houghton, Michigan, 49931 USA\newline
\email{kreher@mtu.edu}
}
\affil[3]{Florida Atlantic University, 
Boca Raton, Florida, 33431 USA\newline
\email{spyros@fau.edu}}
\date{}
\maketitle

\begin{center}{\sl Dedicated to Doug Stinson on the occasion of his 66th birthday.}\end{center}

\abstract{
A \emph{$3$-$(v,\{4,6\},1)$ design} is a configuration 
of  $v$ \emph{points}
and a  collection of $4$- and $6$-element subsets called blocks,
that jointly contain every 3-element subset exactly once.
Using an exhaustive computer search on $v\leq 28$ points
we investigate the $3$-$(v,\{4,6\},1)$ designs that
have a transitive automorphism group and where the
blocks of size 6 form  a 2-class symmetric design.
A 2-class symmetric design with parameters
$(v,k{;}\lambda_1,\lambda_2{;}\delta_1,\delta_2)$
is a set-system on $v$ points and $v$ blocks of size $k$,
where every pair of points
are in $\lambda_1$ or $\lambda_2$ blocks and every pair of blocks intersect in
$\delta_1$ or $\delta_2$ points. The 2-class symmetric designs
include biplanes, semi-biplanes, and 2-class symmetric partially balanced incomplete block designs.
}

\section{Introduction}

A \emph{$3$-$(v,\{4,6\},1)$ design} is a configuration of  $v$ \emph{points}
and a  collection $4$- and $6$-element subsets,
called the \emph{tetrads} and \emph{hexads} respectively,
that jointly contain every 3-element subset exactly once.
Such systems have a long history of investigation, see for example
Hanani~\cite{Hanani},
Mills~\cite{Mills},
Buggenhaut~\cite{Buggenhaut}
and
Kramer-Kreher-Mathon~\cite{KramerKreherMathon}. 
In general a \emph{$t$-wise balanced design} with parameters $t$-$(v,\ensuremath{\mathscr{K}},\lambda)$ is a 
configuration of \emph{points} and subsets of points called \emph{blocks} such that every $t$ element set is contained in $\lambda$ blocks. The parameter  $\ensuremath{\mathscr{K}}$ is the set of allowed block sizes and the parameter $\lambda$ is called the \emph{index}.
If $\ensuremath{\mathscr{K}} = \{k\}$ a $t$-wise balanced design is simply called a $t$-design. 
A $1$-design is also called a \emph{tactical configuration}.
A $t$-design with $b$ blocks and $v$ points is symmetric if $b=v$.
A \emph{$3$-$(v,\{4,6\},1)$ design} 
was termed by Assmus and Sardi~\cite{AssmusSardi} 
to be \emph{homogeneous}
if the tetrads and hexads are each tactical configurations and there are precisely $v$ hexads. 
(Homogeneous is an unfortunate name for such designs,
because it may be confused with homogeneous groups.)
Hence for such designs the hexads form a
symmetric $1$-design with parameters  $1$-$(v,6,6)$.
They prove Theorem~\ref{T-AS}.
\begin{theorem}[Assmus and Sardi~\cite{AssmusSardi}]\label{T-AS}
If a homogenous $3$-$(v,\{4,6\},1)$ design exists, then $v\equiv
2$ or $4 \pmod{6}$ and $v \geq 16$. Moreover, there is a unique such design when $v=16$ and it consists of the best biplane and its sixteen ovals.
\end{theorem}

There are three non-isomorphic $2$-$(16,6,2)$ designs. 
The one with the largest automorphism group was discovered by 
E. Kummer in 1864,~\cite{Kummer}  and has become known as the 
\emph{Kummer configuration}, \cite{Hudson}. It  was studied by 
C. Jordan in 1869, \cite{Jordan} and was called \emph{the best biplane} on 16 points in \cite{AssmusSardi}.

Incidentally, because a homogenous $3$-$(v,\{4,6\},1)$ design contains exactly $v$ hexads, 
it contains $\frac{1}{4}\big(\binom{v}{3}-20v\big)$ tetrads.

The  \emph{automorphism group} $G$ of a $3$-$(v,\{4,6\},1)$ design consists
of those permutations of the points that send blocks to blocks. If this group acts transitively on the points we say that the design is \emph{transitive}.  All designs discussed in this article will be transitive. A transitive
$3$-$(v,\{4,6\},1)$ design with $v$ hexads is automatically homogeneous.
We follow the standard group naming conventions that can for example be found at the url:
\begin{center}
\url{https://people.maths.bris.ac.uk/~matyd/GroupNames/}.
\end{center}
In particular
$\ensuremath{\textsc{A}}_n$, $\ensuremath{\textsc{C}}_n$ and $\ensuremath{\textsc{D}}_n$ are the isomorphism classes of
the alternating, cyclic and dihedral  
groups of degree $n$ and if $q$ is a prime power, then $\ensuremath{\textsc{F}}_q$ denotes a group
isomorphic to the set of all affine 
transformations of  the field $\ensuremath{\mathbb{F}}_q$:
\[
\ensuremath{\textsc{F}}_q = \{ x\mapsto \alpha x + \beta : 
\alpha,\beta \in \ensuremath{\mathbb{F}}_q\text{ and }\alpha\not=0\} \cong \ensuremath{\mathbb{F}}_q \rtimes \ensuremath{\mathbb{F}}_q^{\times}.
\]
It is a Frobenius group.

\section{2-class symmetric designs}

A $1$-$(v,k,k)$ design where each pair of points is in 
exactly $\lambda = \frac{k(k-1)}{v-1}$ blocks is  a 
symmetric design with parameters $(v,k,\lambda)$ or equivalently a
$(v,k,\lambda)$ design or $2$-$(v,k,\lambda)$ design.
A $(v,k,2)$ design is also known as a \emph{biplane} ($\ensuremath{\textsc{bp}}(v,k)$).

A $1$-$(v,k,k)$ design where every pair of points
is either in $\lambda_1$ or $\lambda_2$ blocks 
and every pair of blocks intersect in $\delta_1$ or $\delta_2$ points is 
a $(v,k{;}\lambda_1,\lambda_2{;}\delta_1,\delta_2)$-design. 
The pairs are partitioned into two \emph{classes}.
The pairs occurring in $\lambda_1$ blocks are first class pairs,
whereas the pairs occurring in $\lambda_2$ blocks are 
second class pairs.
We call such designs \emph{2-class symmetric designs}.

A $2$-class symmetric design with $\lambda_1=\lambda_2 = \lambda$ 
is a symmetric design and has $\delta_1=\delta_2 = \lambda$.

A $2$-class symmetric design with $\lambda_1=\delta_1 = 0$ and
$\lambda_2=\delta_2 = 2$ is called a \emph{semi-biplane} ($\ensuremath{\textsc{sbp}}(v,k)$).

A \emph{2-class symmetric design with parameters}
$(v,k{;}\lambda_1,\lambda_2{;}\delta_1,\delta_2)$ is also a
 \emph{2-class partially balanced incomplete block design} 
($\ensuremath{\textsc{pbibd}}(v,k,\lambda_1,\lambda_2)$) if
the two pair classes form a 2-class association scheme.
In this case the first and second class pairs are called 
$i$-th associates for $i=1$ and $2$ respectively.
In such a scheme there are constants $n_1$ and $n_2$
and fixed 2 by 2 matrices $P_1$ and $P_2$. 
Each symbol has exactly $n_i$ $i$-th associates.
If $x$ and $y$ are $i$-th associates, then
the number of other symbols that are $j$-th associates of $x$
and $k$-th associates of $y$ is $P_i[j,k]$.
Only the  $\ensuremath{\textsc{pbibd}}(v,k,\lambda_1,\lambda_2)$s where every pair of blocks intersect in
$\delta_1$ or $\delta_2$ points are $(v,k{;}\lambda_1,\lambda_2{;}\delta_1,\delta_2)$-designs.
Clatworthy in~\cite{Clatworthy} identifies 8 types of 2-class partially balanced incomplete block designs among which only two types are relevant to this article:

\smallskip
\begin{description}
\item[Singular group divisible type:]
These $\ensuremath{\textsc{pbibd}}$s are obtained by starting with a BIBD and then replacing each point by a constant number of new points.
\smallskip
\item[Regular group divisible type:]
A $\ensuremath{\textsc{pbibd}}$ is of this type if $ r-\lambda_2 > 0$ and $rk-v\lambda_1 > 0$,
where $r$ is the number of blocks containing a given point. Note $r=k$ when the design is symmetric.
\end{description}
\smallskip

In this article we enumerate 
transitive $3$-$(v,\{4,6\},1)$ designs, where the hexads form 
a $2$-class symmetric design. These designs will be labeled $D^v_i$
and presented as a list of baseblocks that can be developed by a given group of automorphisms whose generators we provide. (Any listed $3$-$(v,\{4,6\},1)$ design, where the hexads \textbf{do not form} a $2$-class symmetric design are denoted by $X^v_i$.)

\newcommand{\KREHERitem}[1][]{%
  \closepage\item[#1]\minipage[t]{0.95\linewidth}%
  \let\closepage\endminipage%
  }
\newenvironment{KREHERblock-description}{%
  \medskip
  \description
  \let\olditem\item
  \let\closepage\relax  
}{%
  \closepage
  \medskip
  \enddescription
}

\subsection{$v=16$}\label{16}

\begin{minipage}{0.625\textwidth}
The unique homogeneous $3$-$(16,\{4,6\},1)$  design $D_1^{16}$ can be constructed 
from the adjacent 4 by 4 square. For each point $0,1,\ldots,f$ in the square we associate the 6-element blocks consisting of those points in its row and column but not including it. These 
6-element blocks form a 2-(16,6,2) design that is known as the \emph{Kummer configuration} or as the \emph{best biplane on 16 points}. 
\end{minipage}
\hfill
\begin{minipage}{0.18\textwidth}
\begin{center}
\vspace*{-1.25\baselineskip}
\tikz[baseline=(O),scale=0.4,line width=1]{
\coordinate(O) at (1.5,1.5);
\foreach \i in {0,...,3}
{
 \foreach \j in {0,...,3}
 {
 \pgfmathtruncatemacro{\k}{\i + 4 * (3-\j)};
 \coordinate(x\k) at (\i,\j);
 }
}
\node at  (x0) {$0$};
\node at  (x1) {$1$};
\node at  (x2) {$2$};
\node at  (x3) {$3$};
\node at  (x4) {$4$};
\node at  (x5) {$5$};
\node at  (x6) {$6$};
\node at  (x7) {$7$};
\node at  (x8) {$8$};
\node at  (x9) {$9$};
\node at (x10) {$a$};
\node at (x11) {$b$};
\node at (x12) {$c$};
\node at (x13) {$d$};
\node at (x14) {$e$};
\node at (x15) {$f$};
\draw (-0.5,-0.5)-- ++(4,0) -- ++(0,4) -- ++(-4,0) --cycle;
}
\end{center}
\end{minipage}

The 4-element blocks consist of the 4 corners of each of the 36 sub-rectangles of the square and 
the 4 points in each of the 24 transversals. These 60 4-element blocks are the 60 \emph{ovals} of the best biplane.
This exact construction of the 2-(16,6,2) biplane first
appeared in the 1869 paper of C. Jordan where he also determines its automorphism group,~\cite{Jordan}. Assmus and Sardi reveal 
in their 1981 paper~\cite{AssmusSardi} that the blocks 
of size 4 constructed by rectangles and transversals are the biplane's ovals.

\medskip
\begin{description}[font=\upshape\bfseries,leftmargin=2ex]
\item[The best biplane:]\quad\\[0.125\baselineskip]
$12348c$
$02359d$
$0136ae$
$0127bf$
$05678c$
$14679d$
$2457ae$
$3456bf$
$049abc$
$158abd$
$2689be$
$3789af$
$048def$
$159cef$
$26acdf$
$37bcde$
\smallskip
\item[The ovals]\quad\\[-1.125\baselineskip]
\begin{description}[font=\upshape\bfseries,leftmargin=2ex]
\item[From rectangles]
$0145$
$0246$
$0347$
$1256$
$1357$
$2367$
$0189$
$028a$
$038b$
$129a$
$139b$
$23ab$
$01cd$
$02ce$
$03cf$
$12de$
$13df$
$23ef$
$4589$
$468a$
$478b$
$569a$
$579b$
$67ab$
$45cd$
$46ce$
$47cf$
$56de$
$57df$
$67ef$
$89cd$
$8ace$
$8bcf$
$9ade$
$9bdf$
$abef$
\item[From transversals]
$05af$
$05be$
$069f$
$06bd$
$079e$
$07ad$
$14af$
$14be$
$168f$
$16bc$
$178e$
$17ac$
$249f$
$24bd$
$258f$
$25bc$
$278d$
$279c$
$349e$
$34ad$
$358e$
$35ac$
$368d$
$369c$
\end{description}
\smallskip
\end{description}

\begin{KREHERblock-description}
\KREHERitem[$D^{16}_1$:]
\noindent\textbf{Generators}\\
{
(0,1,2,3)(4,5,6,7)(8,9,a,b)(c,d,e,f)\\
(0,4,8,c)(1,5,9,d)(2,6,a,14)(3,7,b,f)\\
}
\noindent\textbf{Baseblocks}\\
{
$\{1,2,3,4,8,c\}$,\\
$\{0,1,4,5\}$,
$\{0,2,4,6\}$,
$\{0,1,8,9\}$,
$\{0,2,8,a\}$,\\
$\{0,7,a,d\}$,
$\{0,5,b,e\}$,
$\{0,5,a,f\}$
}
\end{KREHERblock-description}
The full automorphism group of $D^{16}_1$ is $\ensuremath{\textsc{C}}_2^4.\ensuremath{\textsc{A}}_6.\ensuremath{\textsc{C}}_2$
and is generated by:
\[
(0,2,a)(1,5,3,9,c,7)(6,f,8)(b,e)
\quad\text{and}\quad
(0,f,5,d,c,2)(1,9,e,6,a,7)(3,b,4).
\]

\vbox{
Naturally the $3$-$(16,\{4,6\},1)$ design can also be described bigraphicaly by taking all subgraphs of $K_{4,4}$, isomorphic to
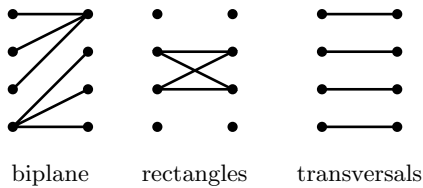
\begin{center}
\begin{tabular}{c@{\qquad}c@{\qquad}c}
\\
\tikz[baseline=(O),scale=1, line width =1]
{ 
 \foreach \i in {0,1,2,3}
 {
  \coordinate(x\i) at (0,0.5*\i); \coordinate(y\i) at (1,0.5*\i); 
  \draw[fill] (x\i) circle[radius=0.05];
  \draw[fill] (y\i) circle[radius=0.05];
 }
  \draw (x0)--(y0) (x0)--(y1) (x0)--(y2) ;
  \draw (y3)--(x1) (y3)--(x2) (y3)--(x3) ;
  \coordinate(O) at ($(x0)!0.5!(y3)$);
}&
\tikz[baseline=(O),scale=1, line width =1]
{ 
 \foreach \i in {0,1,2,3}
 {
  \coordinate(x\i) at (0,0.5*\i); \coordinate(y\i) at (1,0.5*\i); 
  \draw[fill] (x\i) circle[radius=0.05];
  \draw[fill] (y\i) circle[radius=0.05];
 }
 \draw(x1)--(y1)--(x2)--(y2)--(x1);
  \coordinate(O) at ($(x0)!0.5!(y3)$);
}&
\tikz[baseline=(O),scale=1, line width =1]
{ 
 \foreach \i in {0,1,2,3}
 {
  \coordinate(x\i) at (0,0.5*\i); \coordinate(y\i) at (1,0.5*\i); 
  \draw[fill] (x\i) circle[radius=0.05];
  \draw[fill] (y\i) circle[radius=0.05];
  \draw (x\i)--(y\i);
 }
  \coordinate(O) at ($(x0)!0.5!(y3)$);
}\\
\\
biplane&rectangles& transversals\\
\\
\end{tabular}
\end{center}
A complete list of bigraphical $t$-$(v,K,1)$ designs can be found in ~\cite{HoffmanKreher}.
}

Alternatively we can describe the $3$-$(16,\{4,6\},1)$ by
embedding the unique graphical $2$-$(15,\{3,5\},1)$. We take the 
point set to be $E(K_6) \cup \{\infty\}$. 
The blocks are all subgraphs of $K_6$ isomorphic to the
following.

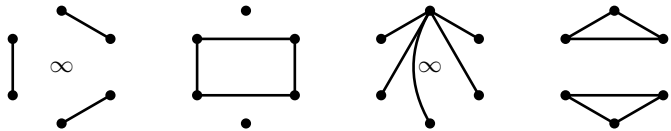
\begin{center}
\begin{tabular}{c*{3}{@{\hspace*{0.4in}}c}}
\\
\tikz[baseline=(O),scale=1, line width =1]
{ 
 \foreach \i in {0,1,...,5}
 {
  \coordinate(x\i) at (90-\i*60:0.75); 
  \draw[fill] (x\i) circle[radius=0.05];
 }
 \draw(x0)--(x1) (x2)--(x3) (x4)--(x5);
  \coordinate(O) at (0,0);
 \node at (O) {$\infty$};
}&
\tikz[baseline=(O),scale=1, line width =1]
{ 
 \foreach \i in {0,1,...,5}
 {
  \coordinate(x\i) at (90-\i*60:0.75); 
  \draw[fill] (x\i) circle[radius=0.05];
 }
 \draw(x1)--(x2)--(x4)--(x5)--(x1);
  \coordinate(O) at (0,0);
}&
\tikz[baseline=(O),scale=1, line width =1]
{ 
 \foreach \i in {0,1,...,5}
 {
  \coordinate(x\i) at (90-\i*60:0.75); 
  \draw[fill] (x\i) circle[radius=0.05];
 }
 \draw(x0)--(x1) (x0)--(x2) (x0)--(x4) (x0)--(x5);
 \draw(x0) to[bend right=30] (x3);
  \coordinate(O) at (0,0);
 \node at (O) {$\infty$};
}&

\tikz[baseline=(O),scale=1, line width =1]
{ 
 \foreach \i in {0,1,...,5}
 {
  \coordinate(x\i) at (90-\i*60:0.75); 
  \draw[fill] (x\i) circle[radius=0.05];
 }
 \draw(x0)--(x1)--(x5)--(x0) (x2)--(x3)--(x4)--(x2);
 \coordinate(O) at (0,0);
}
\end{tabular}
\end{center}
\bigskip

\noindent{}A complete list of graphical $t$-$(v,K,1)$ designs can be found in ~\cite{ChouinardKramerKreher}.

\subsection{$v=20$}\label{20}

Salwach~\cite{Salwach} studied  homogeneous $3$-$(20,\{4,6\},1)$ designs in 1985 and found two whose hexads form a semi-biplane.
His motivation arose by considering a possible extension 
of a supposed projective plane of order 10. Should such an extension (a 3-(112, 12, 1) design) have existed, then it is not difficult to show
that there would necessarily be two blocks whose symmetric difference would consist
of 20 points, see~\cite{Kreher}.
The remaining blocks intersect these 20 points in 2, 4 or 6 points giving
rise to the tetrads and hexads of a 
homogeneous 3-(20, \{4, 6\}, 1), 
see~\cite{AssmusSardi, Buggenhaut} . However in 1989 
it was shown that a projective plane of order 10 cannot 
exist~\cite{Lam}.

Only one of Salwach's two homogeneous 
$3$-$(20,\{4,6\},1)$ designs has a transitive automorphism group.
Salwach reports Andries Brouwer also constructed a
homogeneous $3$-$(20,\{4,6\},1)$ design, where the 
blocks of size 6 do not form a semi-biplane. 
This was a private communication and as far as we can determine 
the Brouwer homogeneous $3$-$(20,\{4,6\},1)$ design was 
never published and we do not know the details of this construction. Andries writes us:
\begin{quotation}
\textit{
I recall the situation with this design. Ed Assmus gave a talk
early in some conference in England, and mentioned that such a
design was unknown. I sat down in the grass a sunny afternoon
and constructed one, and gave him a copy.
}
\end{quotation}
Indeed there is a postscript to the Assmus and Sardi article
\cite{AssmusSardi} that confirms this.
Thus it is likely that the private communication that Salwach
reports was passed on to him from Ed Assmus or refers to this postscript.

It turns out that there are in fact exactly \textbf{three} transitive homogeneous $3$-$(20,\{4,6\},1)$ designs which
we label below as $D^{20}_1$, $X^{20}_2$ and $X^{20}_3$.
Only the hexads of $D^{20}_{1}$,  form a 2-class symmetric 
design.  These hexads are in fact the semi-biplane $\ensuremath{\textsc{sbp}}(20,6)$.
The full automorphism group of these transitive homogeneous $3$-$(20,\{4,6\},1)$ designs is denoted by $G$ in the following table.

\[
\begin{array}{lcllcl}
\hline
&\qquad&G& \vert G\vert &\qquad&\textrm{Hexads type}\\
\hline
D^{20}_{1}&&\bigg.\ensuremath{\textsc{A}}_5\bigg.&60&&\ensuremath{\textsc{sbp}}(20,6)\\
\hline
X^{20}_{2}&&\bigg.\ensuremath{\textsc{F}}_5\bigg.&20&&\textrm{Not a 2-class symmetric design}\\
\hline
X^{20}_{3}&&\bigg.\ensuremath{\textsc{C}}_2^2\times \ensuremath{\textsc{F}}_5\bigg.\qquad&80&&\textrm{Not a 2-class symmetric design}\\
\hline
\end{array}
\]
\begin{KREHERblock-description}
\KREHERitem[$D^{20}_{1}$:]
\noindent\textbf{Generators}\\
{
(0,1)(2,12)(3,13)(4,5)(6,18)(7,19)(8,14)(9,15)(10,16)(11,17)\\
(0,2,8)(1,3,9)(4,11,12)(5,10,13)(6)(7)(14,19,17)(15,18,16)\\
}
\noindent\textbf{Baseblocks}\\
{
$\{0,1,2,14\}$,
$\{0,1,3,11\}$,
$\{0,2,6,8\}$,
$\{0,2,10,13\}$,
$\{0,3,7,17\}$,
$\{0,2,5,7,12,18\}$
}\\
\KREHERitem[$X^{20}_{2}$:]
\noindent\textbf{Generators}\\
{
(0,1)(2,19)(3,18)(4,17)(5,16)(6,14)(7,15)(8,13)(9,12)(10,11)\\
(0,2,8,6)(1,3,9,7)(4,15,5,14)(10,13,18,16)(11,12,19,17)\\
}
\noindent\textbf{Baseblocks}\\
{
$\{0,2,11,14\}$,
$\{0,1,10,11\}$,
$\{0,1,2,8\}$,
$\{0,1,5,7\}$,
$\{0,3,8,11\}$,
$\{0,2,10,12\}$,
$\{0,3,5,16\}$,
$\{0,5,12,15\}$,
$\{0,3,10,18\}$,
$\{0,4,9,16\}$,
$\{0,5,11,19\}$,
$\{0,1,3,4,6,12\}$
}\\
\KREHERitem[$X^{20}_{3}$:]
\noindent\textbf{Generators}\\
{
(0,1)(2,3)(4,5)(6,7)(8,9)(10,11)(12,13)(14,15)(16,17)(18,19)\\
(0,2)(1,3)(4,18)(5,19)(6,16)(7,17)(8,15)(9,14)(10,13)(11,12)\\
(0)(1)(2,7,18,15)(3,6,19,14)(4,13,17,8)(5,12,16,9)(10)(11)\\
}
\noindent\textbf{Baseblocks}\\
{
$\{0,1,10,11\}$,
$\{0,2,5,10\}$,
$\{0,2,7,8\}$,
$\{0,1,4,17\}$,
$\{0,3,11,13\}$,
$\{0,2,11,12\}$,
$\{0,1,2,6,14,18\}$,
}
\end{KREHERblock-description}

\subsection{$v=22$}\label{22}

     There is a unique $\ensuremath{\textsc{sbp}}(22,6)$. It is the double of the 
Paley biplane.  The Paley biplane is the  $2$-$(11,5,2)$
design and is constructed by developing  $\{1,3,4,5,9\}$, 
(the quadratic residues) modulo $11$. If $A$ is an incidence matrix for
the $2$-$(11,5,2)$ then the incidence matrix for the $\ensuremath{\textsc{sbp}}(22,6)$ is
$\left[ \begin{array}{cc}A&I\\I&A^{\textsc{t}}\end{array}\right]$.
This incidence matrix construction is known as \emph{doubling} and is due to Hughes and Dickey, see~\cite{Hughes}.
It can also be obtained by developing 
$\{0_0,1_0,2_0,4_0,7_0,0_1\}$ by the group 
$G$ acting on $\{x_i: x\in\ensuremath{\mathbb{Z}}_{11}, i \in \ensuremath{\mathbb{Z}}_{2} \}$
generated by 
\[
\begin{array}{ccc}
x_i \mapsto (x+1)_i&\text{ and }&x_i \leftrightarrow (-x)_{i+1},
\end{array}
\]
which is isomorphic to the dihedral group $\ensuremath{\textsc{D}}_{11}$,.

There are up to isomorphism 
3 transitive $3$-$(22,\{4,6\},1)$ designs  whose hexads form the
$\ensuremath{\textsc{sbp}}(22,6)$. These are design $D^{22}_{1}$, $D^{22}_{2}$ and $D^{22}_{3}$ 
presented below as a set of 16 
baseblocks with respect to the automorphism group $G$.  Designs 
$X^{22}_{4}$,
$X^{22}_{5}$,
$X^{22}_{6}$, and
$X^{22}_{7}$,
are the transitive $3$-$(22,\{4,6\},1)$ designs that exist on 22 points
whose hexads do not form a 2-class symmetric design. 
They have 19, 20, 20 and 15 baseblocks respectively.

\begin{KREHERblock-description}
\smallskip
\KREHERitem[$D^{22}_{1}$:]\textbf{Baseblocks}\\
{
$\{0_0,2_0,6_0,1_1\}$,
$\{0_0,1_0,2_1,4_1\}$,
$\{0_0,1_0,9_0,3_1\}$,
$\{0_0,3_0,0_1,3_1\}$,
$\{0_0,4_0,2_1,6_1\}$,
$\{0_0,1_0,5_1,6_1\}$,
$\{0_0,1_0,1_1,7_1\}$,
$\{0_0,2_0,2_1,7_1\}$,
$\{0_0,1_0,5_0,8_1\}$,
$\{0_0,3_0,4_1,8_1\}$,
$\{0_0,2_0,6_1,8_1\}$,
$\{0_0,1_0,8_0,9_1\}$,
$\{0_0,4_0,3_1,9_1\}$,
$\{0_0,3_0,6_1,9_1\}$,
$\{0_0,2_0,8_0,10_1\}$,
$\{0_0,1_0,2_0,4_0,7_0,0_1\}$
}\\
\KREHERitem[$D^{22}_{2}$:]\textbf{Baseblocks}\\
{
$\{0_0,1_0,9_0,2_1\}$,
$\{0_0,2_0,6_0,1_1\}$,
$\{0_0,1_0,8_0,3_1\}$,
$\{0_0,3_0,0_1,3_1\}$,
$\{0_0,3_0,1_1,5_1\}$,
$\{0_0,2_0,3_1,5_1\}$,
$\{0_0,1_0,4_1,5_1\}$,
$\{0_0,1_0,1_1,7_1\}$,
$\{0_0,2_0,2_1,7_1\}$,
$\{0_0,1_0,5_0,8_1\}$,
$\{0_0,4_0,1_1,8_1\}$,
$\{0_0,2_0,6_1,8_1\}$,
$\{0_0,4_0,3_1,9_1\}$,
$\{0_0,1_0,6_1,9_1\}$,
$\{0_0,2_0,8_0,10_1\}$,
{$\{0_0,1_0,2_0,4_0,7_0,0_1\}$}
}\\
\KREHERitem[$D^{22}_{3}$:]\textbf{Baseblocks}\\
{
$\{0_0,1_0,9_0,1_1\}$,
$\{0_0,1_0,2_1,7_1\}$,
$\{0_0,1_0,5_0,3_1\}$,
$\{0_0,3_0,1_1,4_1\}$,
$\{0_0,2_0,2_1,4_1\}$,
$\{0_0,1_0,8_0,5_1\}$,
$\{0_0,3_0,0_1,6_1\}$,
$\{0_0,1_0,4_1,6_1\}$,
$\{0_0,2_0,6_0,7_1\}$,
$\{0_0,4_0,3_1,7_1\}$,
$\{0_0,2_0,1_1,8_1\}$,
$\{0_0,3_0,5_1,9_1\}$,
$\{0_0,1_0,8_1,9_1\}$,
$\{0_0,2_0,8_0,10_1\}$,
$\{0_0,5_0,5_1,10_1\}$,
{$\{0_0,1_0,2_0,4_0,7_0,0_1\}$}
}\\
\KREHERitem[$X^{22}_{4}$:]\textbf{Baseblocks}\\
{
$\{0_0,1_0,3_0,5_0\}$,
$\{0_0,1_0,4_0,7_0\}$,
$\{0_0,2_0,7_0,0_1\}$,
$\{0_0,2_0,5_0,1_1\}$,
$\{0_0,1_0,2_1,7_1\}$,
$\{0_0,3_0,1_1,4_1\}$,
$\{0_0,4_0,1_1,5_1\}$,
$\{0_0,1_0,6_0,6_1\}$,
$\{0_0,2_0,4_1,6_1\}$,
$\{0_0,2_0,6_0,8_1\}$,
$\{0_0,5_0,3_1,8_1\}$,
$\{0_0,4_0,2_1,9_1\}$,
$\{0_0,3_0,6_1,9_1\}$,
$\{0_0,2_0,7_1,9_1\}$,
$\{0_0,1_0,2_0,10_1\}$,
$\{0_0,4_0,3_1,10_1\}$,
$\{0_0,5_0,4_1,10_1\}$,
{$\{0_0,1_0,9_0,0_1,1_1,3_1\}$,\,
$\{0_0,1_0,8_0,4_1,5_1,8_1\}$}
}\\
\KREHERitem[$X^{22}_{5}$:]\textbf{Baseblocks}\\
{
$\{0_0,1_0,3_0,6_0\}$,
$\{0_0,1_0,4_0,5_0\}$,
$\{0_0,2_0,4_0,8_0\}$,
$\{0_0,3_0,1_1,4_1\}$,
$\{0_0,5_0,0_1,5_1\}$,
$\{0_0,1_0,4_1,5_1\}$,
$\{0_0,4_0,2_1,6_1\}$,
$\{0_0,2_0,0_1,7_1\}$,
$\{0_0,1_0,2_1,8_1\}$,
$\{0_0,4_0,3_1,7_1\}$,
$\{0_0,1_0,6_1,7_1\}$,
$\{0_0,3_0,0_1,8_1\}$,
$\{0_0,5_0,3_1,8_1\}$,
$\{0_0,4_0,4_1,8_1\}$,
$\{0_0,3_0,6_1,9_1\}$,
$\{0_0,2_0,4_1,9_1\}$,
$\{0_0,1_0,2_0,10_1\}$,
$\{0_0,3_0,7_1,10_1\}$,
{$\{0_0,1_0,9_0,0_1,1_1,3_1\}$,}
{$\{0_0,2_0,7_0,1_1,6_1,8_1\}$}
}\\
\KREHERitem[$X^{22}_{6}$:]\textbf{Baseblocks}\\
{
$\{0_0,1_0,3_0,5_0\}$,
$\{0_0,1_0,4_0,7_0\}$,
$\{0_0,1_0,6_0,8_0\}$,
$\{0_0,3_0,1_1,4_1\}$,
$\{0_0,5_0,0_1,5_1\}$,
$\{0_0,1_0,4_1,5_1\}$,
$\{0_0,4_0,2_1,6_1\}$,
$\{0_0,2_0,0_1,7_1\}$,
$\{0_0,1_0,2_1,8_1\}$,
$\{0_0,4_0,3_1,7_1\}$,
$\{0_0,1_0,6_1,7_1\}$,
$\{0_0,3_0,0_1,8_1\}$,
$\{0_0,5_0,3_1,8_1\}$,
$\{0_0,4_0,4_1,8_1\}$,
$\{0_0,3_0,6_1,9_1\}$,
$\{0_0,2_0,4_1,9_1\}$,
$\{0_0,1_0,2_0,10_1\}$,
$\{0_0,3_0,7_1,10_1\}$,
{$\{0_0,1_0,9_0,0_1,1_1,3_1\}$,},
{$\{0_0,2_0,7_0,1_1,6_1,8_1\}$}
}\\
\KREHERitem[$X^{22}_{7}$:]\textbf{Baseblocks}\\
{
$\{0_0,1_0,4_0,5_0\}$,
$\{0_0,2_0,5_0,7_0\}$,
$\{0_0,1_0,3_0,1_1\}$,
$\{0_0,3_0,7_0,3_1\}$,
$\{0_0,1_0,9_0,3_1\}$,
$\{0_0,1_0,2_0,6_1\}$,
$\{0_0,3_0,0_1,6_1\}$,
$\{0_0,4_0,1_1,6_1\}$,
$\{0_0,3_0,4_1,8_1\}$,
$\{0_0,1_0,7_1,8_1\}$,
$\{0_0,3_0,2_1,9_1\}$,
$\{0_0,3_0,6_0,10_1\}$,
$\{0_0,2_0,8_1,10_1\}$,
$\{0_0,1_0,9_1,10_1\}$,
{$\{0_0,1_0,6_0,0_1,2_1,4_1\}$,}
}\\
\end{KREHERblock-description}
\vspace*{0.5\baselineskip}
The full automorphism group for each design except design $D^{22}_3$ 
is $G\cong \ensuremath{\textsc{D}}_{11}$, and has order $22$.
The full automorphism group for design $D^{22}_3$ has order 110 and is
isomorphic to the group $\ensuremath{\textsc{F}}_{11}$.
It is generated by
\[
\begin{array}{cccc}
x_i \mapsto (x+1)_i,&x_i \leftrightarrow (-x)_{i+1}&\text{ and }&
x_i \mapsto (3x+1)_i
\end{array}
\]

\subsection{$v=26$}\label{26}

On $v=26$ points there are 332 nonisomorphic transitive 
$1$-$(26,6,6)$ designs that do not cover a triple twice.
Although there are seven that are 2-class symmetric designs
 only two of them,  $H_1$ and $H_2$, extend to a 
transitive $3$-$(26,\{4,6\},1)$ design.
\begin{KREHERblock-description}
\KREHERitem[$H_1$] is a  2-class symmetric design with parameters 
$(26,6{;}1,2{;}1,2)$.
\smallskip
\KREHERitem[$H_2$]
 is a 2-class symmetric design with parameters 
$(26,6{;}1,6{;}0,2)$
and also  a singular group divisible type
$\ensuremath{\textsc{pbibd}}(26,6,1,6)$ with parameters 
\[
\begin{array}{llll}
n_1=24, &
n_2=1, &
P_1=\left[\begin{array}{rr} 22& 1\\ 1& 0\end{array}\right], &
P_2=\left[\begin{array}{rr} 24& 0\\ 0&0\end{array}\right]
\end{array}.
\]
This $\ensuremath{\textsc{pbibd}}$  is recoded as S44 on pages 104, 116 and 117 of 
\cite{Clatworthy}.
\end{KREHERblock-description}
There is up to isomorphism a unique $3$-$(26,\{4,6\},1)$ 
design  $D^{26}_1$ that extends $H_1$.
The $\ensuremath{\textsc{pbibd}}$ $H_2$ extends to 6 
nonisomorphic $3$-$(26,\{4,6\},1)$ 
designs $D^{26}_2, D^{26}_3, \ldots, D^{26}_7$. Three have $\ensuremath{\textsc{C}}_{26}$ as their full automorphism group
and three have $\ensuremath{\textsc{C}}_{26} \rtimes \ensuremath{\textsc{C}}_{3}$.
Consequently there are {\large \textbf{7} } nonisomorphic transitive
$3$-$(26,\{4,6\},1)$ designs whose hexads form a 2-class symmetric design.

Develop the following baseblocks modulo 26 to obtain the 4
nonisomorphic  special $3$-$(26,\{4,6\},1)$ designs
that have full automorphism group $\ensuremath{\textsc{C}}_{26}$.

\begin{KREHERblock-description}
\KREHERitem[$D^{26}_{1}$:]\textbf{Baseblocks}\\
{
$\{0,2,7,11\}$,
$\{0,1,4,12\}$,
$\{0,2,6,12\}$,
$\{0,1,3,14\}$,
$\{0,2,5,14\}$,
$\{0,4,9,15\}$,
$\{0,1,10,15\}$,
$\{0,1,7,16\}$,
$\{0,4,11,16\}$,
$\{0,3,13,16\}$,
$\{0,3,9,17\}$,
$\{0,3,6,18\}$,
$\{0,5,13,18\}$,
$\{0,1,8,19\}$,
$\{0,2,10,20\}$,
$\{0,4,13,20\}$,
$\{0,1,18,20\}$,
$\{0,1,2,22\}$,
$\{0,3,10,22\}$,
$\{0,1,5,23\}$,
$\{0,1,17,24\}$,
$\{0,1,6,9,11,13\}$
}\\
\KREHERitem[$D^{26}_{2}$:]\textbf{Baseblocks}\\
{
$\{0,2,4,10\}$,
$\{0,1,6,7\}$,
$\{0,1,8,9\}$,
$\{0,3,7,10\}$,
$\{0,1,2,11\}$,
$\{0,1,12,15\}$,
$\{0,4,11,16\}$,
$\{0,2,6,17\}$,
$\{0,3,9,18\}$,
$\{0,2,12,18\}$,
$\{0,2,9,19\}$,
$\{0,4,12,19\}$,
$\{0,3,6,20\}$,
$\{0,3,11,21\}$,
$\{0,4,9,21\}$,
$\{0,2,16,21\}$,
$\{0,2,14,22\}$,
$\{0,1,16,22\}$,
$\{0,1,5,23\}$,
$\{0,1,3,24\}$,
$\{0,1,4,13,14,17\}$,
$\{0,2,7,13,15,20\}$
}\\
\KREHERitem[$D^{26}_{3}$:]\textbf{Baseblocks}\\
{
$\{0,2,4,10\}$,
$\{0,1,2,9\}$,
$\{0,1,7,11\}$,
$\{0,1,5,12\}$,
$\{0,3,9,14\}$,
$\{0,1,3,15\}$,
$\{0,1,6,16\}$,
$\{0,2,11,17\}$,
$\{0,3,7,18\}$,
$\{0,4,12,18\}$,
$\{0,1,10,19\}$,
$\{0,2,12,19\}$,
$\{0,3,6,20\}$,
$\{0,3,11,21\}$,
$\{0,4,16,21\}$,
$\{0,2,18,21\}$,
$\{0,1,18,22\}$,
$\{0,2,16,23\}$,
$\{0,1,21,23\}$,
$\{0,1,20,24\}$,
$\{0,1,4,13,14,17\}$,
$\{0,2,7,13,15,20\}$
}\\
\KREHERitem[$D^{26}_{4}$:]\textbf{Baseblocks}\\
{
$\{0,1,3,6\}$,
$\{0,2,4,8\}$,
$\{0,1,5,10\}$,
$\{0,1,8,11\}$,
$\{0,1,9,15\}$,
$\{0,4,11,16\}$,
$\{0,3,10,17\}$,
$\{0,2,12,18\}$,
$\{0,1,16,18\}$,
$\{0,4,10,19\}$,
$\{0,2,14,19\}$,
$\{0,1,2,21\}$,
$\{0,3,11,21\}$,
$\{0,4,12,21\}$,
$\{0,3,14,22\}$,
$\{0,2,16,22\}$,
$\{0,1,19,22\}$,
$\{0,1,12,23\}$,
$\{0,2,17,23\}$,
$\{0,1,7,24\}$,
$\{0,1,4,13,14,17\}$,
$\{0,2,7,13,15,20\}$
}
\end{KREHERblock-description}

Develop the following baseblocks modulo 26 to obtain the 3
nonisomorphic  special $3$-$(26,\{4,6\},1)$ designs
that have full automorphism group $\ensuremath{\textsc{C}}_{26}\rtimes\ensuremath{\textsc{C}}_3$.

\begin{KREHERblock-description}
\KREHERitem[$D^{26}_{5}$:]\textbf{Baseblocks}\\
{
$\{0,1,5,9\}$,
$\{0,1,6,11\}$,
$\{0,2,8,12\}$,
$\{0,1,3,15\}$,
$\{0,2,11,16\}$,
$\{0,3,6,17\}$,
$\{0,2,6,18\}$,
$\{0,3,7,18\}$,
$\{0,1,12,18\}$,
$\{0,2,4,19\}$,
$\{0,3,9,19\}$,
$\{0,1,10,19\}$,
$\{0,1,16,20\}$,
$\{0,1,7,21\}$,
$\{0,3,11,21\}$,
$\{0,2,9,21\}$,
$\{0,2,5,22\}$,
$\{0,3,8,22\}$,
$\{0,1,2,23\}$,
$\{0,1,8,24\}$,
$\{0,1,4,13,14,17\}$,
$\{0,2,7,13,15,20\}$
}\\
\KREHERitem[$D^{26}_{6}$:]\textbf{Baseblocks}\\
{
$\{0,1,5,6\}$,
$\{0,1,7,8\}$,
$\{0,1,3,9\}$,
$\{0,2,6,10\}$,
$\{0,3,7,10\}$,
$\{0,1,11,12\}$,
$\{0,3,11,14\}$,
$\{0,4,10,15\}$,
$\{0,4,11,16\}$,
$\{0,2,11,17\}$,
$\{0,1,10,18\}$,
$\{0,2,12,18\}$,
$\{0,2,9,19\}$,
$\{0,4,12,19\}$,
$\{0,3,6,20\}$,
$\{0,4,9,21\}$,
$\{0,2,16,21\}$,
$\{0,2,14,22\}$,
$\{0,2,5,23\}$,
$\{0,1,2,24\}$,
$\{0,1,4,13,14,17\}$,
$\{0,2,7,13,15,20\}$
}\\
\KREHERitem[$D^{26}_{7}$:]\textbf{Baseblocks}\\
{
$\{0,1,7,8\}$,
$\{0,1,3,11\}$,
$\{0,2,8,12\}$,
$\{0,2,9,14\}$,
$\{0,4,10,15\}$,
$\{0,1,12,15\}$,
$\{0,1,6,16\}$,
$\{0,2,11,17\}$,
$\{0,2,6,18\}$,
$\{0,1,10,18\}$,
$\{0,4,12,19\}$,
$\{0,3,10,19\}$,
$\{0,3,6,20\}$,
$\{0,1,9,21\}$,
$\{0,2,16,21\}$,
$\{0,1,5,22\}$,
$\{0,3,18,22\}$,
$\{0,2,19,22\}$,
$\{0,2,5,23\}$,
$\{0,1,2,24\}$,
$\{0,1,4,13,14,17\}$,
$\{0,2,7,13,15,20\}$
}
\end{KREHERblock-description}

\subsection{$v=28$}\label{28}
On $v=28$ points there are 854 nonisomorphic 
transitive $1$-$(28,6,6)$ designs that do not 
cover a triple twice. There are four 
that are 2-class symmetric designs,
but only three $H_1$, $H_2$, $H_3$ 
that extend to  a transitive $3$-$(28,\{4,6\},1)$ design.
\begin{KREHERblock-description}
\smallskip
\KREHERitem[$H_1,H_2$]
are each 2-class symmetric design with parameters $(28,6{;}1,2{;}1,2)$
and are each also  regular group divisible type
$\ensuremath{\textsc{pbibd}}(28,6,1,2)$ with parameters 
\[
\begin{array}{llll}
n_1=24, &
n_2=3, &
P_1=\left[\begin{array}{rr} 20& 3\\ 3& 0\end{array}\right], &
P_2=\left[\begin{array}{rr} 24& 0\\ 2&0\end{array}\right]
\end{array}.
\]
Clatworthy records  these $\ensuremath{\textsc{pbibd}}$ parameters as
R171 on pages 178 and 220 of 
\cite{Clatworthy}, but only one solution is provided.

\KREHERitem[$H_3$]
 is the semi-biplane $\ensuremath{\textsc{sbp}}(28,6)$.
\end{KREHERblock-description}

The Clatworthy $\ensuremath{\textsc{pbibd}}(28,6,1,2)$  recorded as R171 extends to the twelve $3$-$(28,\{4,6\},1)$ designs
$D^{28}_3,D^{28}_4,\ldots,D^{28}_{14}$, the other 
$\ensuremath{\textsc{pbibd}}(28,6,1,2)$ extends to the two $3$-$(28,\{4,6\},1)$ 
designs $D^{28}_1$ and $D^{28}_2$.

The group generators 
\begin{equation}
(0,2,4,\ldots,26)(1,3,5,\ldots,27)
\quad\text{and}\quad
(0,1)(2,3)\cdots(26,27)\label{C2xC14}
\end{equation}
generate $\ensuremath{\textsc{C}}_2 \times \ensuremath{\textsc{C}}_{14}$, which
can be used to develop the baseblocks below 
into the required
$3$-$(28,\{4,6\},1)$ designs
$D^{28}_1,D^{28}_2,\ldots, D^{28}_{10}$.
\begin{KREHERblock-description}
\KREHERitem[$D^{28}_{ 1}$:]\textbf{Baseblocks}\\
$\{0,1,4,5\}$,
$\{0,2,8,11\}$,
$\{0,1,10,12\}$,
$\{0,4,9,13\}$,
$\{0,3,4,14\}$,
$\{0,2,13,14\}$,
$\{0,3,9,15\}$,
$\{0,1,14,15\}$,
$\{0,2,9,16\}$,
$\{0,3,11,18\}$,
$\{0,4,10,19\}$,
$\{0,4,8,20\}$,
$\{0,5,14,20\}$,
$\{0,7,15,20\}$,
$\{0,5,8,21\}$,
$\{0,2,12,22\}$,
$\{0,4,17,22\}$,
$\{0,2,4,23\}$,
$\{0,3,12,23\}$,
$\{0,5,17,23\}$,
$\{0,1,18,23\}$,
$\{0,2,7,24\}$,
$\{0,3,8,24\}$,
$\{0,2,15,25\}$,
$\{0,1,8,26\}$,
$\{0,2,5,27\}$,
$\{0,1,2,6,17,20\}$
\KREHERitem[$D^{28}_{ 2}$:]\textbf{Baseblocks}\\
$\{0,1,4,5\}$,
$\{0,2,7,10\}$,
$\{0,3,6,12\}$,
$\{0,2,8,12\}$,
$\{0,2,5,16\}$,
$\{0,3,10,15\}$,
$\{0,2,11,15\}$,
$\{0,1,14,15\}$,
$\{0,3,8,17\}$,
$\{0,4,11,17\}$,
$\{0,2,14,18\}$,
$\{0,1,8,19\}$,
$\{0,3,9,19\}$,
$\{0,5,10,20\}$,
$\{0,7,12,20\}$,
$\{0,2,9,21\}$,
$\{0,6,15,21\}$,
$\{0,1,12,22\}$,
$\{0,3,18,22\}$,
$\{0,3,4,23\}$,
$\{0,2,4,24\}$,
$\{0,3,11,25\}$,
$\{0,4,12,25\}$,
$\{0,2,19,25\}$,
$\{0,1,10,27\}$,
$\{0,2,23,27\}$,
$\{0,1,2,6,17,20\}$
\KREHERitem[$D^{28}_{ 3}$:]\textbf{Baseblocks}\\
$\{0,1,4,5\}$,
$\{0,2,6,11\}$,
$\{0,1,10,11\}$,
$\{0,3,8,12\}$,
$\{0,2,9,12\}$,
$\{0,1,14,15\}$,
$\{0,4,10,16\}$,
$\{0,5,11,16\}$,
$\{0,3,15,16\}$,
$\{0,5,10,17\}$,
$\{0,2,15,17\}$,
$\{0,2,8,18\}$,
$\{0,3,11,18\}$,
$\{0,4,15,18\}$,
$\{0,3,10,19\}$,
$\{0,5,12,20\}$,
$\{0,6,15,20\}$,
$\{0,2,4,21\}$,
$\{0,1,8,22\}$,
$\{0,2,13,22\}$,
$\{0,4,17,22\}$,
$\{0,2,5,25\}$,
$\{0,3,20,24\}$,
$\{0,2,23,24\}$,
$\{0,3,4,25\}$,
$\{0,1,12,26\}$,
$\{0,2,10,27\}$,
$\{0,1,2,7,16,20\}$
\KREHERitem[$D^{28}_{ 4}$:]\textbf{Baseblocks}\\
$\{0,1,4,5\}$,
$\{0,2,8,11\}$,
$\{0,1,10,11\}$,
$\{0,3,7,12\}$,
$\{0,2,9,12\}$,
$\{0,2,5,13\}$,
$\{0,3,9,13\}$,
$\{0,1,14,15\}$,
$\{0,3,15,16\}$,
$\{0,4,10,17\}$,
$\{0,5,11,17\}$,
$\{0,2,15,17\}$,
$\{0,3,8,18\}$,
$\{0,2,6,18\}$,
$\{0,4,15,18\}$,
$\{0,3,10,19\}$,
$\{0,4,9,20\}$,
$\{0,5,12,20\}$,
$\{0,6,15,20\}$,
$\{0,2,10,21\}$,
$\{0,1,8,22\}$,
$\{0,5,10,22\}$,
$\{0,3,4,24\}$,
$\{0,2,19,24\}$,
$\{0,2,4,25\}$,
$\{0,1,12,26\}$,
$\{0,2,22,27\}$,
$\{0,1,2,7,16,20\}$
\KREHERitem[$D^{28}_{ 5}$:]\textbf{Baseblocks}\\
$\{0,1,4,5\}$,
$\{0,3,4,10\}$,
$\{0,2,5,10\}$,
$\{0,1,10,11\}$,
$\{0,2,8,12\}$,
$\{0,1,14,15\}$,
$\{0,3,15,16\}$,
$\{0,4,9,17\}$,
$\{0,2,15,17\}$,
$\{0,2,6,18\}$,
$\{0,3,11,18\}$,
$\{0,4,15,18\}$,
$\{0,5,8,20\}$,
$\{0,6,15,20\}$,
$\{0,2,9,21\}$,
$\{0,5,11,21\}$,
$\{0,1,8,22\}$,
$\{0,3,7,22\}$,
$\{0,2,13,22\}$,
$\{0,3,13,23\}$,
$\{0,5,16,23\}$,
$\{0,2,19,24\}$,
$\{0,3,20,24\}$,
$\{0,2,4,25\}$,
$\{0,3,19,25\}$,
$\{0,1,12,26\}$,
$\{0,2,11,27\}$,
$\{0,1,2,7,16,20\}$
\KREHERitem[$D^{28}_{ 6}$:]\textbf{Baseblocks}\\
$\{0,1,4,5\}$,
$\{0,2,6,11\}$,
$\{0,1,10,11\}$,
$\{0,3,8,12\}$,
$\{0,2,9,12\}$,
$\{0,1,14,15\}$,
$\{0,4,10,16\}$,
$\{0,5,11,16\}$,
$\{0,5,10,17\}$,
$\{0,2,14,17\}$,
$\{0,2,8,18\}$,
$\{0,3,11,18\}$,
$\{0,4,15,18\}$,
$\{0,3,10,19\}$,
$\{0,5,12,20\}$,
$\{0,6,15,20\}$,
$\{0,2,4,21\}$,
$\{0,1,8,22\}$,
$\{0,2,13,22\}$,
$\{0,4,17,22\}$,
$\{0,2,5,25\}$,
$\{0,3,20,24\}$,
$\{0,2,23,24\}$,
$\{0,3,4,25\}$,
$\{0,2,10,27\}$,
$\{0,1,12,27\}$,
$\{0,1,2,7,16,20\}$
\KREHERitem[$D^{28}_{ 7}$:]\textbf{Baseblocks}\\
$\{0,1,4,5\}$,
$\{0,2,8,11\}$,
$\{0,1,10,11\}$,
$\{0,3,7,12\}$,
$\{0,2,9,12\}$,
$\{0,2,5,13\}$,
$\{0,3,9,13\}$,
$\{0,1,14,15\}$,
$\{0,4,10,17\}$,
$\{0,5,11,17\}$,
$\{0,2,14,17\}$,
$\{0,3,8,18\}$,
$\{0,2,6,18\}$,
$\{0,4,15,18\}$,
$\{0,3,10,19\}$,
$\{0,4,9,20\}$,
$\{0,5,12,20\}$,
$\{0,6,15,20\}$,
$\{0,2,10,21\}$,
$\{0,1,8,22\}$,
$\{0,5,10,22\}$,
$\{0,3,4,24\}$,
$\{0,2,19,24\}$,
$\{0,2,4,25\}$,
$\{0,1,12,27\}$,
$\{0,2,22,27\}$,
$\{0,1,2,7,16,20\}$
\KREHERitem[$D^{28}_{ 8}$:]\textbf{Baseblocks}\\
$\{0,1,4,5\}$,
$\{0,3,4,10\}$,
$\{0,2,5,10\}$,
$\{0,1,10,11\}$,
$\{0,2,8,12\}$,
$\{0,1,14,15\}$,
$\{0,4,9,17\}$,
$\{0,2,14,17\}$,
$\{0,2,6,18\}$,
$\{0,3,11,18\}$,
$\{0,4,15,18\}$,
$\{0,5,8,20\}$,
$\{0,6,15,20\}$,
$\{0,2,9,21\}$,
$\{0,5,11,21\}$,
$\{0,1,8,22\}$,
$\{0,3,7,22\}$,
$\{0,2,13,22\}$,
$\{0,3,13,23\}$,
$\{0,5,16,23\}$,
$\{0,2,19,24\}$,
$\{0,3,20,24\}$,
$\{0,2,4,25\}$,
$\{0,3,19,25\}$,
$\{0,2,11,27\}$,
$\{0,1,12,27\}$,
$\{0,1,2,7,16,20\}$
\KREHERitem[$D^{28}_{ 9}$:]\textbf{Baseblocks}\\
$\{0,1,4,5\}$,
$\{0,2,6,11\}$,
$\{0,1,10,11\}$,
$\{0,3,8,12\}$,
$\{0,2,9,12\}$,
$\{0,1,14,15\}$,
$\{0,4,10,16\}$,
$\{0,5,11,16\}$,
$\{0,3,15,16\}$,
$\{0,5,10,17\}$,
$\{0,2,15,17\}$,
$\{0,2,8,18\}$,
$\{0,3,11,18\}$,
$\{0,4,15,18\}$,
$\{0,3,10,19\}$,
$\{0,5,12,20\}$,
$\{0,6,14,20\}$,
$\{0,2,4,21\}$,
$\{0,7,14,21\}$,
$\{0,2,13,22\}$,
$\{0,4,17,22\}$,
$\{0,1,8,23\}$,
$\{0,2,5,25\}$,
$\{0,3,20,24\}$,
$\{0,2,23,24\}$,
$\{0,3,4,25\}$,
$\{0,1,12,26\}$,
$\{0,2,10,27\}$,
$\{0,1,2,7,16,20\}$
\KREHERitem[$D^{28}_{10}$:]\textbf{Baseblocks}\\
$\{0,2,6,11\}$,
$\{0,3,8,12\}$,
$\{0,2,9,12\}$,
$\{0,1,14,15\}$,
$\{0,4,10,16\}$,
$\{0,5,11,16\}$,
$\{0,5,10,17\}$,
$\{0,2,14,17\}$,
$\{0,1,4,18\}$,
$\{0,2,8,18\}$,
$\{0,3,11,18\}$,
$\{0,3,10,19\}$,
$\{0,5,12,20\}$,
$\{0,6,15,20\}$,
$\{0,2,4,21\}$,
$\{0,1,8,22\}$,
$\{0,2,13,22\}$,
$\{0,4,17,22\}$,
$\{0,2,5,25\}$,
$\{0,3,20,24\}$,
$\{0,2,23,24\}$,
$\{0,3,4,25\}$,
$\{0,1,10,25\}$,
$\{0,2,10,27\}$,
$\{0,1,12,27\}$,
$\{0,1,2,7,16,20\}$
\end{KREHERblock-description}

The remaining 4  transitive $3$-$(28,\{4,6\},1)$ 
are obtained by developing the baseblocks below with
the group $\ensuremath{\textsc{F}}_8$ generated by
\begin{equation}
{\footnotesize
\left.
\begin{array}{@{}l@{}}
(0,4,8,13,17,21,25)(1,5,9,12,16,20,24)(2,6,10,15,19,23,26)(3,7,11,14,18,22,27)\\
\text{and }(2,3)(4,18)(5,19)(6,7)(8,23)(9,22)(10,25)(11,24)(12,26)(13,27)(16,17)(20,21)
\end{array}\right\}
}
\label{F8}
\end{equation}

\begin{KREHERblock-description}
\KREHERitem[$D^{28}_{11}$:]\textbf{Baseblocks}\\
$\{0,1,4,5\}$,
$\{0,2,4,10\}$,
$\{0,2,6,11\}$,
$\{0,5,8,17\}$,
$\{0,1,14,15\}$,
$\{0,4,11,17\}$,
$\{0,2,14,17\}$,
$\{0,2,13,18\}$,
$\{0,4,15,18\}$,
$\{0,5,15,19\}$,
$\{0,6,15,20\}$,
$\{0,2,9,21\}$,
$\{0,1,8,22\}$,
$\{0,2,5,24\}$,
$\{0,2,19,25\}$,
$\{0,1,12,27\}$,
$\{0,2,8,27\}$,
$\{0,1,2,7,16,20\}$
\KREHERitem[$D^{28}_{12}$:]\textbf{Baseblocks}\\
$\{0,1,4,5\}$,
$\{0,2,4,10\}$,
$\{0,2,6,11\}$,
$\{0,5,8,17\}$,
$\{0,1,14,15\}$,
$\{0,4,11,17\}$,
$\{0,2,15,17\}$,
$\{0,2,13,18\}$,
$\{0,4,15,18\}$,
$\{0,5,15,19\}$,
$\{0,6,15,20\}$,
$\{0,2,9,21\}$,
$\{0,1,8,22\}$,
$\{0,2,5,24\}$,
$\{0,2,19,25\}$,
$\{0,1,12,26\}$,
$\{0,2,8,27\}$,
$\{0,1,13,27\}$,
$\{0,1,2,7,16,20\}$
\KREHERitem[$D^{28}_{13}$:]\textbf{Baseblocks}\\
$\{0,2,4,10\}$,
$\{0,2,6,11\}$,
$\{0,5,8,17\}$,
$\{0,1,14,15\}$,
$\{0,4,11,17\}$,
$\{0,2,14,17\}$,
$\{0,1,4,18\}$,
$\{0,2,13,18\}$,
$\{0,6,15,20\}$,
$\{0,2,9,21\}$,
$\{0,1,8,22\}$,
$\{0,2,5,24\}$,
$\{0,1,10,25\}$,
$\{0,2,19,25\}$,
$\{0,1,12,27\}$,
$\{0,2,8,27\}$,
$\{0,1,2,7,16,20\}$
\KREHERitem[$D^{28}_{14}$:]\textbf{Baseblocks}\\
$\{0,1,4,5\}$,
$\{0,2,4,10\}$,
$\{0,2,6,11\}$,
$\{0,5,8,17\}$,
$\{0,1,14,15\}$,
$\{0,4,11,17\}$,
$\{0,2,15,17\}$,
$\{0,2,13,18\}$,
$\{0,4,15,18\}$,
$\{0,5,15,19\}$,
$\{0,6,14,20\}$,
$\{0,2,9,21\}$,
$\{0,1,9,22\}$,
$\{0,1,8,23\}$,
$\{0,2,5,24\}$,
$\{0,2,19,25\}$,
$\{0,1,12,26\}$,
$\{0,2,8,27\}$,
$\{0,1,13,27\}$,
$\{0,1,2,7,16,20\}$
\end{KREHERblock-description}
The full automorphism group of 
$D^{28}_1,D^{28}_2,\ldots, D^{28}_{10}$ 
has order 28 and is $\ensuremath{\textsc{C}}_2 \times \ensuremath{\textsc{C}}_{14}$ with 
generators~(\ref{C2xC14}).

The full automorphism group of $D^{28}_{11}$ and $D^{28}_{12}$ has order order 56
and is $\ensuremath{\textsc{F}}_8$ with generators~(\ref{F8}), but the full automorphism group of $D^{28}_{13}$ and $D^{28}_{14}$ has order order 168.
It is $\ensuremath{\textsc{F}}_8\rtimes \ensuremath{\textsc{C}}_3 = 
\ensuremath{\textsc{Aut}}(\ensuremath{\textsc{F}}_8)$
with generators~(\ref{F8}) and (\ref{gamma}).
\begin{equation}
{\footnotesize
\left.
\begin{array}{@{}l@{}}
(1,14,15)(2,4,22)(3,18,9)(5,8,17)(6,26,25)(7,12,10)(11,21,27)(13,24,20)(16,19,23)
\end{array}\right.
}
\label{gamma}
\end{equation}
Enumerating  the extensions of the semi-biplane $\ensuremath{\textsc{sbp}}(28,6)$ 
$H_3$ 
has proven to be very difficult, because there are apparently millions
of nonisomorphic extensions. 
We provide baseblocks and generators
for one of the transitive $3$-$(28,\{4,6\},1)$ designs, whose hexads are $H_3$,
the semi-biplane $\ensuremath{\textsc{sbp}}(28,6)$.

\begin{KREHERblock-description}
\KREHERitem[$D^{28}_{15}$:]
\noindent\textbf{Generators}\\
{
\footnotesize
$(4,8,17)(5,9,16)(6,10,18)(7,11,19)(12,25,20)(13,24,21)(14,27,22)(15,26,23)$\\
$(0,1)(2,3)(4,5)(6,7)(8,9)(10,11)(12,13)(14,15)\cdots(24,25)(26,27)$\\
$(0,2,1,3)(4,15,9,27,17,23,5,14,8,26,16,22)(6,13,11,25,18,21,7,12,10,24,19,20)$\\
$(0,4,8,12,17,20,25)(1,5,9,13,16,21,24)(2,6,10,15,18,23,26)(3,7,11,14,19,22,27)$
}
\noindent\textbf{Baseblocks}\\
{
$\{0,2,4,22\}$,
$\{0,1,14,15\}$,
$\{0,1,4,5\}$,
$\{0,2,5,15\}$,
$\{0,4,14,19\}$,
$\{0,4,9,21\}$,
$\{0,2,7,24\}$,
$\{0,5,18,22\}$,
$\{0,4,12,26\}$,
$\{0,4,13,27\}$,
$\{0,1,2,6,10,18\}$
}
\end{KREHERblock-description}

\section{Computation}
Let $G$ be a transitive subgroup of $\textsc{Sym}(X)$. 
If $x \in X$ and $g \in G$, then the image of $x$ under $g$ is
$x^g$ and if $S \subseteq X$, then $S^g = \{x^g : x\in S\}$ 
is the image of $S$. The orbit of $S$ under the action of $G$
is $S^G = \{S^g : \text{$g \in G$}\}$.
Let
\[
\begin{array}{c@{,\ldots,}cl} 
\Delta_1&\Delta_{m}&\text{ be the orbits of $3$-element sets}\\
\Gamma_1&\Gamma_{n}&\text{ be the orbits of $4$-element sets}\\
E_1&E_{s}&\text{ be the orbits of $6$-element sets}
\end{array}
\]
Choose $E_{i_1},E_{i_2},\ldots,E_{i_x}$ such that
$\sum_{j=1}^x \vert E_{i_j}\vert = v$, then 
$\ensuremath{\mathscr{H}}=\cup_{j=1}^x E_{i_j} $ is a $1$-$(v,6,6)$ design and a possible
candidate set of hexads. 
Choose $T_i \in \Delta_i$ to be any fixed representative, $i=1,2,..,m$
and let
\begin{align*}
\ensuremath{\mathscr{R}} &=\{\Delta_i : \vert\{S \in \ensuremath{\mathscr{H}}: S\supset T_i\}\vert =0\},\\
\ensuremath{\mathscr{C}}&=\{\Gamma_j : \vert\{K \in \Gamma_j : K \supset T_i\}\vert \leq 1 \text{ for every }
\Delta_i \in \ensuremath{\mathscr{R}} \}
\end{align*}
Define the $(0,1)$-valued $\ensuremath{\mathscr{R}}$ by $\ensuremath{\mathscr{C}}$ matrix $A$ by
\[
A[\Delta,\Gamma] = \vert \{ K \in \Gamma : K \supset T \}\vert,
\]
where $T \in \Delta$ is a fixed representative. Then
a (0,1) solution $U$ to $AU=J$, where $J=[1,1,\ldots,1]^T$,
selects the tetrads that together with the candidate hexads 
forms a homogeneous $3$-$(v,\{4,6\},1)$ on $v=\vert X\vert$ points that 
has $G$ as an automorphism group. The orbit representatives and construction
of the incidence matrix used tools described in~\cite{KreherStinson} with suitable modification.  The matrix equation  $AU=J$ was solved with a
backtrack algorithm using the dancing links data structure described in
\cite{Knuth}. The implementation we employed can be downloaded from

\begin{center}
\url{http://pottonen.kapsi.fi/libexact.html}
\end{center}

A subgroup $G$ of $\textsc{Sym}(X)$ is \emph{minimally transitive}
if it contains no proper transitive subgroup. To find all transitive designs
on a given number of points $v$ only the minimal groups need be searched.
The following magma code finds generators for all minimal transitive groups
on $v$ points:\\

\begin{center}
\small
\begin{alltt}
procedure MinTrans(v)
 fn:="MinTransGroupsDegree"*IntegerToString(v)*".txt";
 SetAutoColumns(false);
 SetColumns(0);
 for i in [1..NumberOfTransitiveGroups(v)] do
  G,c:=TransitiveGroup(v,i);
  max:=MaximalSubgroups(G);
  if not &or [IsTransitive(max[j]\`{}subgroup): j in [1..#max]] then
   Write(fn, c);
   Write(fn, Sprint(#Generators(G))*" " *Sprint(v));
   for j in [1..NumberOfGenerators(G)] do
    Write(fn,G.j);
   end for;
   Write(fn, GroupName(G));
   Write(fn, "------------");
  end if;
 end for;
 SetAutoColumns(true);
end procedure;
\end{alltt}
\quad\\
\end{center}

For each $v$  we collect all designs  afforded 
by each of the minimal transitive groups of degree $v$.
These designs are then 
reduced to a list of nonisomorphic designs using the graph isomorphism program \verb+nauty+ obtained  here:

\begin{center}
\url{http://users.cecs.anu.edu.au/~bdm/}
\end{center}
Once all designs are found, \verb+nauty+ was used to find their 
full automorphism groups. Then essentially the same method described above was re-run to obtain the baseblocks for the designs with respect to their full groups.

\section*{Acknowledgments}
We thank  the
Fields Institute for Research in Mathematical Sciences
for the support they provided and Rosemary Bailey for 
her very useful comments.

\end{document}